\documentclass[11]{siamltex}
\usepackage{amsmath,amssymb,amsfonts,graphicx,fancybox,shadow,color}
\newtheorem{remark}[theorem]{Remark}

\usepackage{graphicx}

\def\be#1\ee{\begin{equation}#1\end{equation}}

\begin{document}
\bibliographystyle{plain}

\pagestyle{myheadings}
%\markboth{V.~Druskin, V.~Simoncini and M.~Zaslavsky}{Adaptive rational Krylov space}

\title{Reduced order modeling inversion of monostatic data in a multi-scattering environment. }
%\title{LIPPMANN-SCHWINGER-LANCZOS ALGORITHM FOR INVERSE
%	SCATTERING WITH NON-RECIPROCALLY EXTENDED MATRIX TRANSFER FUNCTIONS}
\author{
 V. Druskin\footnotemark[1],  S. Moskow\footnotemark[2] and M. Zaslavsky\footnotemark[3]\textsuperscript{ ,}\footnotemark[4]}

\renewcommand{\thefootnote}{\fnsymbol{footnote}}

\footnotetext[1]{Worcester Polytechnic Institute, Department of Mathematical Sciences,
Stratton Hall,
100 Institute Road, Worcester MA, 01609 (vdruskin1@gmail.com)}
\footnotetext[2]{Department of Mathematics, Drexel University, Korman Center, 3141 Chestnut Street, Philadelphia, PA 19104
(moskow@math.drexel.edu)}
\footnotetext[3]{Schlumberger-Doll Research Center, 1 Hampshire St., 
Cambridge, MA 02139-1578 (mzaslavsky@slb.com)}
\footnotetext[4]{Southern Methodist University, Department of Mathematics, Clements Hall, 3100 Dyer st,
Dallas, TX 75205 (mzaslavskiy@smu.edu)}
\maketitle

\begin{abstract} The data-driven reduced order models (ROMs) have recently emerged as an efficient tool for the solution of the inverse scattering problems with applications to seismic and sonar imaging.  One specification of this approach is that it requires the full square  multiple-output/multiple-input (MIMO) matrix valued transfer function as data for multidimensional problems.  The synthetic aperture radar (SAR), however, is limited to single input/single output (SISO) measurements corresponding to the diagonal of the matrix transfer function.  Here we present  a ROM based   Lippmann-Schwinger  approach overcoming this  drawback. The ROMs are constructed to match the data for each source-receiver pair separately, and these are used to construct internal solutions for the corresponding source using only the data-driven Gramian.  Efficiency of the proposed approach is demonstrated on 2D and 2.5D (3D propagation and 2D reflectors) numerical examples. The new algorithm   not only suppresses multiple echoes seen in the Born imaging, but also takes advantage of  illumination by them of some back sides of the reflectors,  improving the quality of their mapping.
	
 \end{abstract}
\section{Introduction}
The reduced order model (ROM) approach has been shown previously to be a powerful tool for inverse impedance, scattering and diffusion \cite{borcea2011resistor,borcea2014model,druskin2016direct,druskin2018nonlinear,borcea2017untangling,borcea2019robust,BoDrMaMoZa,
borcea2020reduced, DrMoZa,  Borcea2021ReducedOM, DrMoZa2, BoGaMaZi}.  In this work, we apply it via the Lippmann-Schwinger-Lanczos algorithm  \cite{DrMoZa}  to models of synthetic aperture radar (SAR). In the process, we present a simplification of the ROM approach that applies for general time domain problems.  This simplification both makes the algorithm  and its exposition more direct, while yielding potential for computational speedup.

In the ROM framework for solving multidimensional inverse problems, given a full symmetric matrix transfer function, the ROM is chosen precisely to match the given data set, see
\cite{druskin2016direct,druskin2018nonlinear,borcea2017untangling,borcea2019robust, BoDrMaMoZa, Borcea2021ReducedOM, borcea2021reduced}.  Then, the ROM is transformed to a sparse form (tridiagonal for single input single output (SISO) problems, block tridiagonal for multiple input/output (MIMO) problems) by Lanczos orthogonalization. 
This data-driven ROM, in this orthogonalized form, has entries for which their dependence on the unknown PDE coefficients is approximately linear \cite{borcea2011resistor,druskin2016direct,borcea2014model,borcea2020reduced,Borcea2021ReducedOM}.
This process is related to works of Marchenko, Gelfand, Levitan and Krein on inverse spectral problems and to the idea of spectrally matched second order staggered finite-difference grids first introduced in \cite{druskin1999gaussian}, and first used for direct inversion in \cite{BoDr}.  

The data-driven ROM can be viewed as  Galerkin matrix \cite{doi:10.1137/1.9781611974829.ch7}. 
The crucial step, first noticed in \cite{druskin2016direct}, is that the orthogonalized Galerkin basis depends very weakly on the unknown medium.  In \cite{BoDrMaMoZa}, it was shown that this basis allows one to generate internal solutions from boundary data only, and in \cite{DrMoZa}, the data-generated internal solutions ${\bf u}_p$  (corresponding to unknown coefficient $p$) were used in the Lippmann-Schwinger integral equation, yielding the so-called Lippmann-Schwinger-Lanczos  (LSL) method.   In \cite{BoGaMaZi}, this method was first used in the time domain. 

The LSL method for the frequency-domain problems works as follows. Given data $F_p$ corresponding to unknown $p$ and background data $F_0$,  the Lippmann-Schwinger integral equation says that \be\label{eq:LipSwi} F_p -F_0 =-\langle u_0,p u_p\rangle \ee
where $\langle ,\rangle$  is an appropriate inner product on the domain, and where $u_p$ and $u_0$ are the unknown and background internal solutions respectively. We then use the data generated internal solution ${\bf u}_p$ in place of $u_p$:
\be\label{eq:LipSwiL} F_p -F_0  \approx -\langle u_0,p {\bf u}_p\rangle.\ee
Recall that ${\bf u}_p$ is precomputed directly from the data without knowing $p$.

\begin{figure}
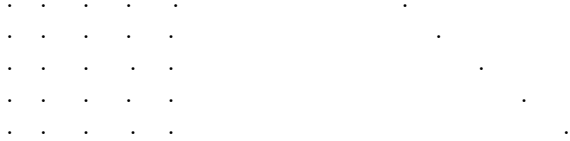

\be\label{eq:ext.matrix} 
	 \begin{matrix}. & .\ &. \ & .\ &.  \\
  .  & . \ &.\  &. \  & .\  \\
     .  &.  \  & .  \  &.   &.  \  \\
   .   & . \  & . \ & .\  &.  \  \\
 .    & . \  & . \  & .  & . \  \\ 
	\end{matrix}  \hskip 1in \ \ \ \begin{matrix}.& \ & \ & \ &  \\
    & . \ &\  & \  & \  \\
       &  \  & .  \  &   &  \  \\
      &  \  &  \ & .\  &  \  \\
     &  \  &  \  &   & . \  \\ 
	\end{matrix} \nonumber  \ee
\label{dataarray}
\caption{Data arrays structures: Full MIMO square transfer function (data from all sources at all receivers)  (left) vs monostatic data (data received only back at each source, typical of SARs) (right).} 
\end{figure}

Although a ROM must be constructed from a full symmetric transfer function, the Lippmann-Schwinger-Lanczos approach allows for its application to more general data sets. In \cite{DrMoZa2}, we showed several examples with non-symmetric data arrays. Here we consider the monostatic formulation, that is, the case where one has only the diagonal of the matrix transfer function. We are again capitalizing on the fact that the Lippmann-Schwinger-Lanczos algorithm does not necessarily require the full matrix transfer function \cite{DrMoZa2}. The structure of our monostatic measurement array can be summarized with the help of Figure \ref{dataarray}.  The column and row numbers correspond respectively to the indices of the receivers (outputs) and transmitters (inputs). The conventional data-driven ROM requires data from all receivers for all transmitters (the full data set), that is, that the measurements are given by full square matrix (Figure \ref{dataarray}).  We assume here that we have only the diagonal part. We construct a separate ROM for each transmitter to obtain its corresponding internal field. Subsequently, all data will be coupled via the Lippmann-Schwinger equation. We show that even with this sparse data array, we retain some of the good performance of the ROM approach for strong nonlinear scattering, where the Born approximation fails.

Somewhat independently of the extension of LSL to monostatic problems, we present a simplification of the generation of the internal solutions in the time domain. In the spectral domain version \cite{borcea2020reduced,DrMoZa,DrMoZa2},  one would compute the mass and stiffness matrix from the data, execute the Lanczos algorithm to orthogonalize, then use the ROM as a forward operator to generate internal solutions.  In the time domain, it turns out that we can obtain the same internal solutions directly from the mass matrix, with no need to introduce finite difference time stepping or the stiffness matrix. The formula we obtain offers potential to speed up the computation, and allows us to present the process more compactly.

This paper is organized as follows. In Section 2 we describe the entire time domain process in detail for a one dimensional single input single output (SISO) problem including numerical experiments.   A brief discussion of this simplified process for  MIMO arrays in higher dimensions is described in Section 3.  Section 4 contains a detailed exposition of the monostatic formulation, and Section 5 contains numerical experiments.

\section{One dimensional problem} 
\subsection{Statement of inverse problem}
Consider the following one-dimensional model problem 
\begin{eqnarray} \label{wavemodel1}
u_{tt}  + Au&= &  0 \ \mbox{in}  \ \Omega\times [ 0,\infty) \\
u ( t=0) &=& g \ \mbox{in}\  \Omega  \\
u_t ( t=0) &=&  0\  \mbox{in}\  \Omega \end{eqnarray}
where $\Omega = (0,1)$, operator \begin{equation}\label{Aeq}  A = A_0 + q \end{equation}
where  $A_0\ge 0$  is known, (for example $A_0=-\Delta $),  $q(x)\ge 0$ is our unknown potential, and initial data is given by 
\begin{equation} \label{geq}
g(x) = \sqrt{ \hat{f}(\sqrt{A})} \delta_0(x) 
\end{equation}
for $\hat{f}(\omega)$ the Fourier transform of the initial pulse, which in our case we choose to be modulated Gaussian
\begin{equation} \label{feq}
f(t) =e^{-\sigma^2t^2/2}\cos(\omega_0t).
\end{equation}
We also assume homogeneous Neumann boundary conditions on the spatial boundary $\partial \Omega$. 
The exact forward solution to (\ref{wavemodel1}) is 
\begin{equation}\label{exactsolution}  u(x,t) = \cos{(\sqrt{A}t)} g(x) . \end{equation}
Of course this is unknown for our inverse problem, except near the receiver.  We measure data at $x=0$ at the $2n-1$ evenly spaced time steps $t= k\tau$ for $k=0,\dots, 2n-2$ integrated against $g$
\begin{equation}\label{datadef} F(k\tau ) = \int_\Omega g(x) \cos{(\sqrt{A}k\tau)} g(x) dx .\end{equation}
Recall that $g$ is concentrated near $x=0$, representing both source and receiver. The inverse problem is as follows:  Given $$ \{ F(k \tau) \} \  \mbox{ for } \  k=0,\dots, 2n-2,$$
reconstruct $q$. 
\begin{remark} \label{initialdataremark} The initial data $g$ is chosen in the form (\ref{geq}) for convenience; without the outer square root in the definition of $g$, the problem  (\ref{wavemodel1}) is equivalent to having zero initial data and source term $\delta_0(x) f(t)$. Furthermore, for our choice of pulse $f$, we have \be\label{eq:modgauss} \hat{f}(\omega) = \frac{\sqrt{\pi}}{\sqrt{2}\sigma}\left(e^{-\frac{(\omega-\omega_0)^2}{2\sigma^2}}+e^{-\frac{(\omega+\omega_0)^2}{2\sigma^2}}\right),\ee
so that for $\omega_0=0$ \begin{equation} \label{ourg} g=\left(\frac{\sqrt{2\pi}}{\sigma}\right)^{1/2} e^{-\frac{A}{2\sigma^2}} \delta_0(x)\end{equation} is the solution of a time domain diffusion equation with operator $A$ at time $\frac{1}{2\sigma^2}$, which is assumed to be early enough that our operator $A$ is still equal to our known $A_0$. The square root in the definition of $g$ is chosen for symmetry reasons which we will see below. \end{remark}
\subsection{Mass matrix}
Define $u_k$ to be the true snapshot $u_k = u(k\tau,x)$  for $k=0,\ldots,2n-2$. Then the $n\times n$  mass matrix is defined by
\begin{equation} \label{massmatrix} M_{kl}= \int_\Omega u_k u_l  dx \end{equation}
for $k,l =0,\ldots,n-1$.  From our expression (\ref{exactsolution}) for the exact solution,
\begin{equation} \label{massmatrix2} M_{kl} = \int_\Omega  g(x) \cos{(\sqrt{A}k\tau)} \cos{(\sqrt{A}l\tau)} g(x) dx. \end{equation}
The cosine angle sum formula 
$$ \cos{(\sqrt{A}k\tau)} \cos{(\sqrt{A}l\tau)} = \frac{1}{2} \left(  \cos{(\sqrt{A}(k-l) \tau)} + \cos{(\sqrt{A}(k+l) \tau)} \right) $$
along with (\ref{datadef}) gives us directly that 
\begin{equation} \label{massfromdata} M_{kl} = \frac{1}{2} \left( F((k-l)\tau) + F((k+l)\tau) \right),\end{equation} 
that is, we can obtain this mass matrix directly from the data. Note that we need precisely the $2n-1$ data points corresponding to $k=0,\dots, 2n-2$ to obtain this $n\times n$ mass matrix. 
\subsection{Orthogonalization}  Note that $M$ is positive definite, so we can compute its Cholesky decomposition 
$$M=U^\top U$$
where $U$ is upper triangular. 
This would give us an orthogonalization of the true snapshots, if we knew them, as follows. Define $\vec{u}$ to be a row vector of the first $n$ snapshots corresponding to $k=0,\ldots,n-1$, and set 
\begin{equation}\label{vdef}  \vec{v} = \vec{u} U^{-1} , \end{equation}
that is,  we set \begin{equation} \label{vdef2} v_k = \sum_l u_l U^{-1}_{lk}. \end{equation}
Then the set of functions $\{ v_k \}$ will be orthonormal in the $L^2$ norm. One can check
\begin{eqnarray} \int_\Omega v_i v_j dx   &=& \int_\Omega \left( \sum_l U^{-1}_{li} u_l  \right)\left( \sum_k U^{-1}_{kj} u_k  \right) dx \nonumber \\ &=& \sum_{lk} U^{-1}_{li} U^{-1}_{kj}  \int_\Omega u_l u_k dx \nonumber \\
&=& \sum_{lk} U^{-1}_{li} U^{-1}_{kj}  M_{lk} \nonumber \\
&=& \sum_{lkr} U^{-1}_{li} U^{-1}_{kj}  U^\top_{lr}U_{rk} \nonumber \\
&=& \delta_{ij}. \nonumber \end{eqnarray}
Since $U^{-1}$ is upper triangular, this procedure is precisely Gram-Schmidt performed on the true snapshots, in order of the time steps. 
That is, although we do not know the true snapshots, the mass matrix, obtained directly from the data, gives the transformation that orthogonalizes them.  
\subsection{Background problem} We now do all of the above for the known background problem,  which has exact solution \begin{equation}\label{exactsolution0}  u^0(x,t) = \cos{(\sqrt{A_0}t)} g(x). \end{equation}
We have the corresponding background snapshots $\{ u_j^0 \} $, mass matrix
\begin{equation} \label{massmatrix0} M^0_{kl}= \int_\Omega u^0_k u^0_l  dx, \end{equation}
corresponding Cholesky decomposition $$M^0=(U^0)^\top U^0,$$
and orthogonalized background snapshots
\begin{equation}\label{v0def}  \vec{v}^0 = \vec{u}^0 (U^0)^{-1} . \end{equation}
\subsection{Crucial step} 
{ It was noticed in  \cite{druskin2016direct} that  the {\it orthogonalized} snapshots depend very weakly on $q$.  That is
\begin{equation}\label{eq:as.independence}  \vec{v} \approx \vec{v}^0.\end{equation}
The idea why is that since we are doing Gram-Schmidt \eqref{vdef2} on the time snapshots in sequential order, we are orthogonalizing away any reflections, since they are overlapping in space with previous times. A rigorous analysis for a related problem  for a so-called optimal grid is given in \cite{borcea2005continuum}.  The Gram-Schmidt procedure \eqref{vdef2} is closely related to the Marchenko-Gelfand-Levitan setting,  to be examined more fully later.}    

\subsection{Data-driven internal solutions}  
\begin{figure}[htb]
	\centering
	\begin{tabular}{cc}
		\includegraphics[width=0.45\textwidth]{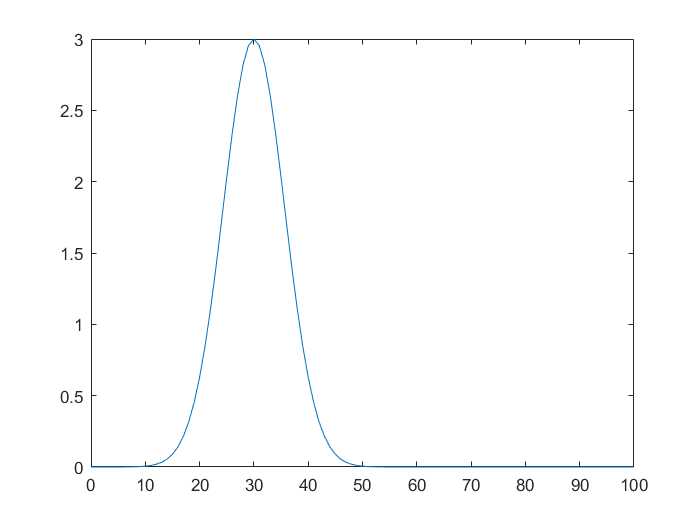} &
		\includegraphics[width=0.45\textwidth]{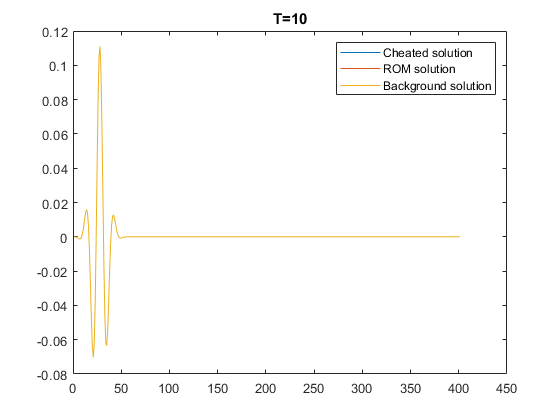} \\
		\includegraphics[width=0.45\textwidth]{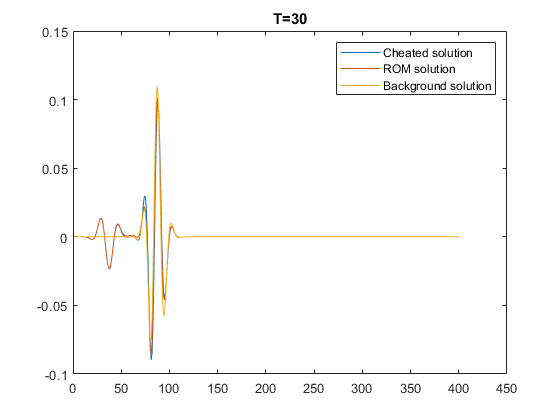}&
		\includegraphics[width=0.45\textwidth]{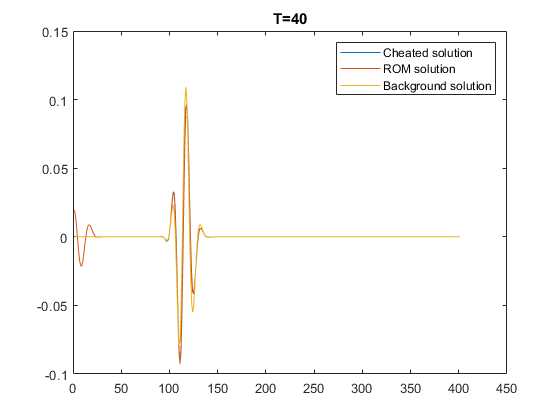} \\
	\end{tabular}
	\caption{Data generated internal snapshots compared to the true snapshots and the background snapshots. Before hitting the scatterer, all overlap. After hitting the scatterer, the data generated solutions show reflections very close to the true  ("cheated")  solutions, while the background solutions are quite different.}
	\label{fig:internal1d1}
\end{figure}

From (\ref{vdef}) and (\ref{eq:as.independence}) we have that the true snapshots 
\begin{eqnarray} \vec{u}  &=& \vec{v} U\nonumber \\
     &\approx&  \vec{v}^0 U. \nonumber\end{eqnarray}
This motivates the definition of our data generated snapshots 
  \begin{eqnarray}  \vec{\bf{u}} &=&  \vec{v}^0U \nonumber \\ &=&\vec{u}^0 (U^0)^{-1}U. \label{internal} \end{eqnarray}
Note that in the right hand side above we needed only the background solutions and the true $U$, which we obtained just from data. Equation (\ref{internal}) is a simple formula for the data generated internal solutions. 
In Fig. \ref{fig:internal1d1} we illustrated how accurately $\vec{\bf{u}}$ approximates true snapshots $\vec{u}$ and how different both are compared to the background snapshots $\vec{u}^0$. Here we probed a medium with one bump (see top left in Fig. \ref{fig:internal1d1} by modulated Gaussan waveform (\ref{eq:modgauss}) with $\frac{\omega_0}{\sqrt{\sigma}}=3$ excited and measured at $x=0$. As one can expect, all three snapshots are on top of each other for small times when the waveform didn't reach the reflector yet (see top right in Fig. \ref{fig:internal1d1}). In turn, for later times the background snapshots become totally different from the true snapshots, however the latter are still reproduced by ROM-generated snapshots $\vec{\bf{u}}$ quite accurately (see bottom left and bottom right in Fig. \ref{fig:internal1d1}).

\subsection{Inversion using Lippmann-Schwinger-Lanczos}
For each time step $k\tau$, $k= 0,\dots, n-1$, consider the time domain Lippmann-Schwinger equation  (see for example \cite{BoGaMaZi})
\begin{equation} \label{Lipsch} F_0(k\tau) - F(k\tau) = \int_0^{k\tau}\int_\Omega u_0(x,k\tau-t) u(x,t) q(x) dxdt .\end{equation}
Recall for the Born approximation one would replace $u(x,t)$ with the background $u_0(x,t)$ in the integral above. We instead replace $u$ by its (time semi-discretized) data generated approximation $\vec{\bf{u}}$ from (\ref{internal}). Let ${\bf u}(x,t)$ be a time interpolant of $\vec{\bf{u}}$ , to produce the Lippmann-Schwinger-Lanczos equation
\begin{equation} \label{LipschL} F_0(k\tau) - F(k\tau) = \int_0^{k\tau} \int_\Omega u_0(x,k\tau-t) {\bf{u}}(x,t) q(x) dxdt ,\end{equation}
which can be inverted for $q$. 
We note that instead of interpolating the snapshots in time, equivalently one may use the discrete time steps as nodes in some numerical approximation of (\ref{Lipsch}).
\subsection{Remarks on the relation to previous work}
 In \cite{druskin2016direct}, it was shown that snapshot $u_i$ can be obtained by the action of the $i$-th (first kind) Chebyshev polynomial of the propagation operator $P=\cos(\tau \sqrt A)$ on the initial condition $g$.  Hence the snapshots $u_i$ form nested Krylov subspaces, that is  $$ \hbox{span}\{u_0, \ldots, u_{i}\} = \hbox{span}\{g, Pg,\ldots, P^{i} g\},$$ for all $i=0,\ldots n-1$. Therefore, Gram-Schmidt orthogonalization on the snapshots is equivalent to the Lanczos algorithm,  hence the name Lippmann-Schwinger-Lanczos (LSL).  The LSL algorithm was first introduced  in \cite{DrMoZa}, which was in the frequency domain, in which case the Lanczos algorithm is executed explicitly. 

 If one wants to compute a full ROM, we need a stiffness matrix too. The reconstruction methods in much of the previous work used the full ROM \cite{druskin2016direct,druskin2018nonlinear,borcea2019robust,borcea2021reduced}. We see here, though, that for time domain Lippmann-Schwinger-Lanczos, we need only the internal solutions, and hence only the mass matrix $M$ is needed, yielding the above algorithm simplification.  There is also potential for speedup due to fast decay of the row elements of the upper triangular matrix $(U^0)^{-1}U$, potentially allowing for truncation. This is to be investigated in future work.  

\section{Full multidimensional problem-MIMO formulation}
All of the above simplification extends to MIMO problems in higher dimensions. Consider 
\begin{eqnarray} \label{wavemodel2}
u_{tt}  + Au&= &  0 \ \mbox{in}  \ \Omega\times [ 0,\infty) \\
u ( t=0) &=& g \ \mbox{in}\  \Omega  \\
u_t ( t=0) &=&  0\  \mbox{in}\  \Omega \end{eqnarray}
where here $\Omega $ is a domain in $\mathbb{R}^d$, and again operator  $A = A_0 + q$. Consider the a set of source/receivers modeled by $\{ g_j \}$, a pulse localized near a receiver point $x_j$, and the response matrix
\begin{equation}\label{datadef2} F^{ji} (k\tau) = \int_\Omega g_j(x) \cos{(\sqrt{A}k\tau)} g_i(x) dx ,\end{equation}
for $j,i=1,\dots, m$, $k=0,\ldots, 2n-1$, representing the response at receiver $j$ from source $i$ at time $k\tau$. 
For the full MIMO problem, the mass tensor can again be obtained by the extension of (\ref{massfromdata}) to blocks (as in many earlier works, see for example \cite{BoGaMaZi}), and the internal solutions again obtained from the block Cholesky  (with $m\times m$ blocks) decomposition, 
$$ M = U^\top U$$ 
where $U$ is block upper triangular.  We note that there is some ambiguity in the choice of the blocks; we choose them so that all resulting orthogonalized functions 
$$ \vec{v} = \vec{u} U^{-1} .$$
are all mutually orthogonal. We do a similar decomposition/orthogonalization for the background mass matrix
$$ M^0 = (U^0)^\top U^0, \ \ \ \  \vec{v}^0 = \vec{u}^0 (U^0)^{-1} $$
to obtain the data generated internal solutions directly 
$$\vec{\bf{u}} = \vec{u}^0 (U^0)^{-1} U .$$

\section{Monostatic formulation} 
\subsection{Linear-algebraic setup}
 For the synthetic aperture radar, we are given only the diagonal response 
\begin{equation}\label{datadefradar} F^{jj} (k\tau) = \int_\Omega g_j(x) \cos{(\sqrt{A}k\tau)} g_j(x) dx ,\end{equation}
for $j=1,\dots, m$, $k=0,\ldots, 2n-1$, from which we are not able to obtain the complete mass matrix.  However, we can instead compute a mass matrix and corresponding internal solution corresponding to each source separately.

Let $u^{(j)}$ be the true solution given data $g_j$. Define its set of snapshots $$u^{(j)}_k = u^{(j)}(k\tau,x)$$  for $k=0,\ldots,2n-1$, $j=1,\dots,m$ to correspond to source $j$ at time $k\tau $. For SARs this is read only at receiver $j$. Then the $n\times n$  mass matrix corresponding to this source is defined by
\begin{equation} \label{massmatrixj} M^j_{kl}= \int_\Omega u^j_k u^j_l  dx \end{equation}
for $k,l =0,\ldots,n-1$. Again we have the expression (\ref{exactsolution}) for the exact solution, and
\begin{equation} \label{massmatrixj2} M^j_{kl} = \int_\Omega  g_j(x) \cos{(\sqrt{A}k\tau)} \cos{(\sqrt{A}l\tau)} g_j(x) dx, \end{equation}
so the cosine angle formula 
and (\ref{datadefradar}) gives us directly that 
\begin{equation} \label{massfromdataj} M^{(j)}_{kl} = \frac{1}{2} \left( F^{jj}((k-l)\tau) + F^{jj}((k+l)\tau) \right).\end{equation} 
Each $M^{(j)}$ is positive definite, so we can compute its Cholesky decomposition 
$$M=(U^{(j)})^\top{U^{(j)}}$$
where each $U^{(j)}$ is upper triangular, for $j=1,\ldots, m$. 
The orthogonalized true snapshots corresponding to source $j$ are obtained in the exact same way,
\begin{equation}\label{vdefj}  \vec{v}^{(j)} = \vec{u}^{(j)}( U^{(j)})^{-1} , \end{equation}
where $ \vec{u}^{(j)}$ is the column vector of the original (yet unknown) snapshots.
So we again have that \begin{equation} \label{vdef2j} v^{(j)}_k = \sum_l u^{(j)}_l  {(U^{(j)})}^{-1}_{lk} , \end{equation}
and the set of functions $\{ v^{(j)}_k \}$ for $k= 0,1,\dots, n-1$ will be orthonormal in the $L^2$ norm. 
We now do all of the above for the known background problem,  which has exact solutions
\begin{equation}\label{exactsolution0j}  u^{(j),0}(x,t) = \cos{(\sqrt{A_0}t)} g(x). \end{equation}
We have the corresponding background snapshots $\{ u^{(j),0}\} $, mass matrices
\begin{equation} \label{massmatrix0j} M^{(j),0}_{kl}= \int_\Omega  u^{(j),0}_k  u^{(j),0}_l dx, \end{equation}
corresponding Cholesky decompositions $$M^{(j),0}=(U^{(j),0})^\top U^{(j),0},$$
and orthogonalized background snapshots
\begin{equation}\label{v0defj}  \vec{v}^{(j),0} = \vec{u}^{(j),0} (U^{(j),0})^{-1} . \end{equation}
\subsection{Crucial step}
 To accurately approximate the internal solution , we need an analogue of \eqref{eq:as.independence}. It was shown in \cite{druskin2016direct, druskin2018nonlinear} that \eqref{eq:as.independence} still holds for the MIMO formulation, provided we have sufficient array density and aperture. In contrast, the monostatic setup does not provide enough functions in the subspace to cancel all of the reflections from different directions during orthogonalization. However, it will still cancel the most important reflections,  by the following reasoning.  
 
Consider the  $3D$ problem in the half-space, i.e.,  with $\Omega=(\mathbb{R}^3)^+$. Then for regular enough $q$ and sharp pulse (small $\sigma$), we have that for a given source,
 \begin{equation} \label{eq:anzac}
 u(x,t)\approx \frac{1}{2\pi}\frac{\delta(\|x\|-t)}{\|x\|}+\eth(x,t),\end{equation}
 where  $\eth(x,t)$ is a smooth  function satisfying the causality principle $\eth=0$ if $\|x\|\ge t$ \cite{}. Then the Gram-Schmidt orthogonalization of  $u_i =u(x,i\tau)$ to  $u_j=u(x,j\tau)$ for $j<i$  approximately cancels the spherical averages of $\eth(x,t)$.
 Thus, the  {\it spherical averages} of orthogonalized snapshots depend very weakly on $q$.  That is, for each $j=1,\dots,m$,
\begin{equation} \label{spherical} \int_{|x|=const}\left[\vec{v}^{(j)} -\vec{v}^{(j),0}\right]\approx 0. \end{equation}

\subsection{Data generated internal solutions}
The formulas (\ref{vdefj}), (\ref{v0defj}) and (\ref{spherical}) yields our data generated internal snapshots for each $j$,
  \begin{eqnarray} \label{internalj} \vec{\bf{u}}^{(j)} &=&\vec{v}^{(j),0} U^{(j)}  \\ &=& \vec{u}^{(j),0} (U^{(j),0})^{-1}U^{(j)} . \nonumber\end{eqnarray}
  These will approximate the true solution in the sense of of averages on spheres centered at the transmitter/receiver location. 
  
 \subsection{Monostatic Lippmann-Schwinger-Lanczos equation}
 We now consider the Lippmann-Schwinger equation. For each time step $k\tau$, $k= 0,\dots, n-1$, and for each source $j=1,\ldots, m$, we have 
\begin{equation} \label{Lipschj} F^{jj}_{0}(k\tau) - F^{jj}(k\tau) = \int_0^{k\tau}\int_\Omega u^{(j),0}(x,k\tau-t) u^{(j)}(x,t) q(x) dxdt .\end{equation}
For the Lippmann-Schwinger-Lanczos method, we replace $u^{(j)}$ by its data generated approximation $\vec{\bf{u}}^{(j)}$ from (\ref{internalj}), which we again need to interpolate. Let ${\bf u}^{(j)}(x,t)$ be a time interpolant of $\vec{\bf{u}}^{(j)}$ , to produce the Lippmann-Schwinger-Lanczos equation
\begin{equation} \label{LipschLj} F^{jj}_0(k\tau) - F^{jj}(k\tau) = \int_0^{k\tau}\int_\Omega u^{(j),0}(x,k\tau-t) {\bf{u}}^{(j)}(x,t) q(x) dxdt ,\end{equation}
for $j=1,\ldots,m$ , $k= 0,\ldots ,n-1$, yielding $nm$ equations to be inverted for $q$. 
The spheres discussed in the context of approximation \eqref{internalj} coincide with 
with  the slowness surfaces from the transmitter/receiver location. Thus the corresponding spherical waves give the dominant contribution in \eqref{Lipschj},and so it is reasonable to expect that \eqref{LipschLj} is a good approximation of \eqref{Lipschj}. Below we will verify this reasoning via numerical experiments.

\section{Numerical experiments}

\begin{figure}
	\centering
	\begin{tabular}{cc}
		\includegraphics[width=0.46\textwidth]{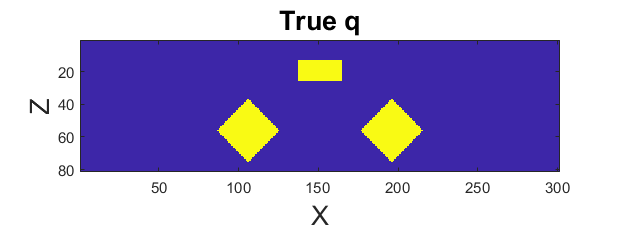} &
		\includegraphics[width=0.46\textwidth]{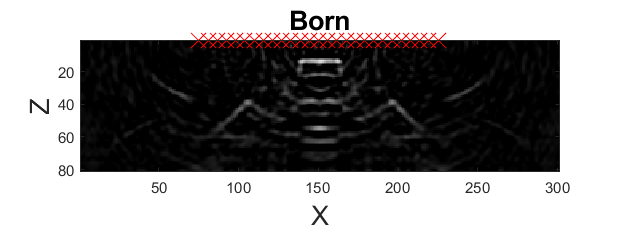} 	\end{tabular}
		\includegraphics[width=0.46\textwidth]{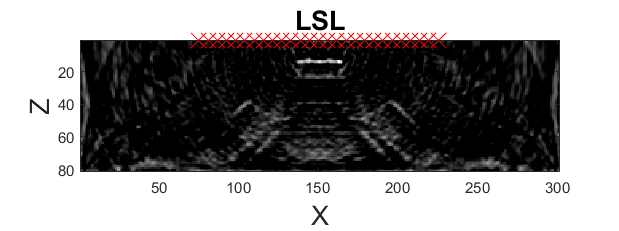}

	\caption {Numerical experiment 1:  2D inverse scattering problem  with 3 massive reflectors.  True medium (top left) and its reconstructions using the Born linearization (top right), and our LSL(bottom). The  monostatic array locations are shown in red. We image the reflectivity distribution to map the  boundaries of the reflectors. The LSL algorithm    not only suppresses multiple echoes seen in the Born imaging but also takes advantage of  illumination by them of some back sides of the reflectors  improving quality of their mapping.  The noise seen at the periphery   of the LSL image  will be addressed in the future by improved processing.}
	\label{fig:numex2}
\end{figure}

We start with the  2D inverse scattering problem to image 3 infinitely conductive reflectors shown on the top right in Fig.~\ref{fig:numex2}. We note that in this multi-scattering environment it is highly challenging to reconstruct all of the details using SAR data only. We generate synthetic data by discretizing (\ref{wavemodel2}) in $\Omega=[0;300]\times [0;80]$ on a $600\times 160$ grid and then solving the obtained equations for 27 sources emitting a non-modulated Gaussian pulse. Lippmann-Schwinger equation was approximated using quadrature on $300\times 80$ grid. The Born solution captures the top boundaries of all three reflectors nicely, however the bottom boundaries remain invisible. Also,  the top reflector produced multiple ghost images (see the top right plot in Fig.~\ref{fig:numex2}). In turn, LSL managed not only to avoid that, but also was able to exploit multiple scattering effects to map most of the bottom boundaries, as well as the internal structure of the large reflectors (see the bottom plot in Fig. \ref{fig:numex2}). 
 
\begin{figure}[htb]
	\centering
		\begin{tabular}{cc}
		\includegraphics[width=0.42\textwidth]{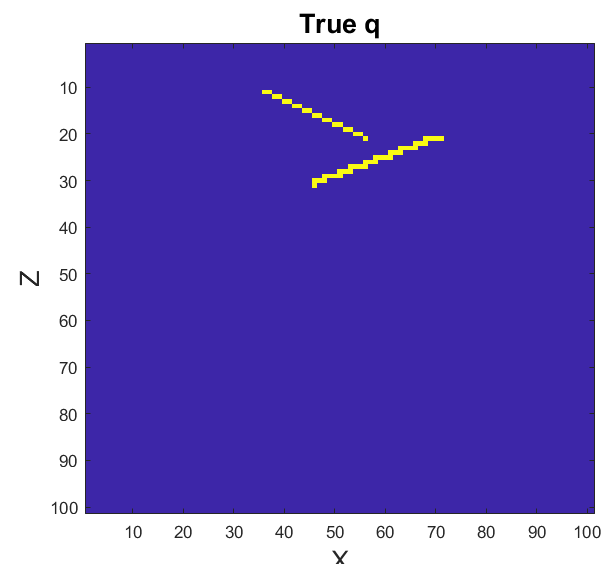} &
		\includegraphics[width=0.42\textwidth]{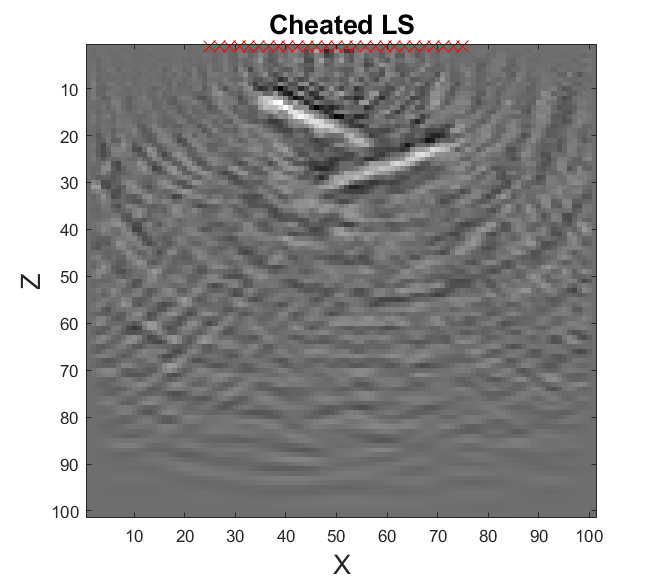}\\
		\includegraphics[width=0.42\textwidth]{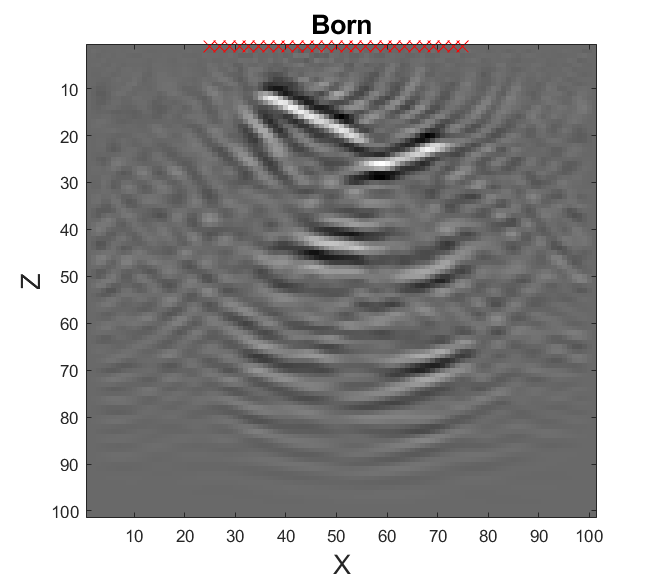} &
		\includegraphics[width=0.42\textwidth]{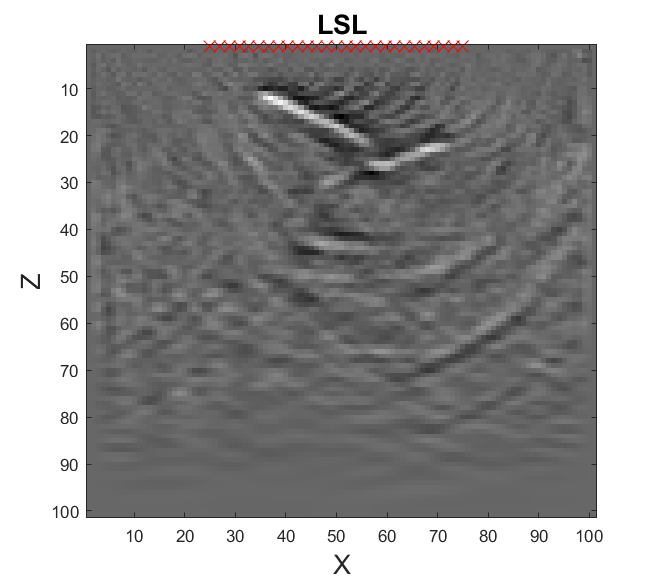}\\
	\end{tabular}
	\caption{{Numerical experiment 2:  2.5D inverse scattering problem  with 2 elongated thin reflectors. True medium (top left) and its reconstructions using the cheated IE (top right), Born linearization (bottom left) and our LSL (bottom right). Similar to the experiment 1, the  monostatic array locations are shown in red.}}
	\label{fig:numex1}
\end{figure}

In the second numerical example we consider so-called 2.5D SAR inverse scattering problem, that is, a 3D wave  scattered by 2D  reflectors.  Here, we consider  2 thin elongated reflectors  embedded in a homogeneous background (see the top left in Fig. \ref{fig:numex1}). For the SAR data collected along a single trajectory, it is a reasonable approximation to assume that the medium is uniform in a horizontal direction perpendicular to the trajectory. To obtain synthetic data, we discretize (\ref{wavemodel2}) in $\Omega=[0;100]^3$ using finite-differences on a $200\times 200\times 200$ 3D grid. The obtained discrete problem was solved for 27 positions of radar that were emitting a non-modulated Gaussian pulse (radar position are marked by crosses in Fig. \ref{fig:numex1}). Then the Lippmann-Schwinger equation (\ref{LipschLj}) is approximated using quadrature on a $100\times 100$ 2D grid. On the top right in Fig. \ref{fig:numex1} we plot a cheated version of the Lippmann-Schwinger solution, meaning that all internal solutions $u^{(j)}(x,t)$ are assumed to be  known exactly. We note that this corresponds to the best case scenario that one could expect from solving (\ref{Lipschj}) iteratively, that is, by updating the background solution computed from the approximate $q$ obtained in the previous step.\footnote{ Such an Iterative Lippman-Schwinger (a.k.a. distorted Born) method) would require multiple solutions of the forward problems, which can be prohibitively  expensive for radar imaging even if the iterations converge.}.  In the bottom two plots of Fig. \ref{fig:numex1} we have shown the Born solution and the solution produced by our LSL approach respectively. Due to the lack of aperture in SAR data, both approaches failed to image a part of the lower reflector. However,  the Born solution also suffers from multiple ghost images of reflectors that are caused by the misinterpretation of multi-scattering effects. In turn, the LSL solution managed to clean out multiple echoes and produced a significantly better image. In fact, this image is just slightly inferior to the one from the "cheated IE", which indicates that we have a good approximation of the interior solution.  

\section{Conclusion}
In this work we describe how the LSL (Lippmann-Schwinger-Lanczos) approach simplifies in the time domain, and demonstrate its usefulness for sparse data sets, in particular for the monostatic (SAR) problem. We exploit the fact the Lippmann-Schwinger equation allows us to plug in the an approximate internal solutions computed separately from each SISO data set, and then use all of the given data simultaneously in the LSL system. We compared the solution produced with our LSL approach to the Born solution for examples in 2 and 2.5 dimensions. Due to the lack of aperture in SAR data, both approaches fail to image all of the reflector.  However, the LSL approach eliminates much of the ghost images in the Born solution, in addition to imaging somewhat more of the reflectors. 

$$ $$
\thanks{{\bf Acknowledgements} 
This material is partially based upon work supported by the NSF grant DMS-1929284 while authors were in residence at the ICERM in Providence, RI, during the Spring 2020 Semester Program "Model and dimension reduction in uncertain and dynamic systems" and Spring 2020 Reunion Event. Authors are also grateful to Liliana Borcea, Alex Mamonov and J\"orn Zimmerling for productive discussions that inspired this research. V. Druskin was partially supported by AFOSR grant FA 955020-1-0079 {and NSF grant  DMS-2110773}. S. Moskow was partially supported by NSF grant DMS-2008441.  }

\bibliography{biblio,biblio6,graphbib6,galerkincitations}
\end{document}